\documentclass[11pt]{article}
\usepackage{amssymb,amsfonts,amsmath,latexsym,epsf,tikz,url}

\newtheorem{theorem}{Theorem}[section]

\newtheorem{conjecture}[theorem]{Conjecture}
\newtheorem{corollary}[theorem]{Corollary}

\newtheorem{definition}[theorem]{Definition}

\newcommand{\proof}{\noindent{\bf Proof.\ }}
\newcommand{\qed}{\hfill $\square$\medskip}

\begin{document}

\title{An upper bound on the distinguishing index of graphs with minimum degree at least two}

\author{Saeid Alikhani   $^{}$\footnote{Corresponding author} \and Samaneh Soltani
}
\date{\today}

\maketitle

\begin{center}
Department of Mathematics, Yazd University, 89195-741, Yazd, Iran\\
{\tt  alikhani@yazd.ac.ir, s.soltani1979@gmail.com}
\end{center}


\begin{abstract}
The distinguishing index of a simple graph $G$, denoted by $D'(G)$, is the
least number of labels in an edge labeling of $G$ not preserved by
any non-trivial automorphism.  It was conjectured by Pil\'sniak (2015)
that for any 2-connected graph $D'(G) \leq \lceil  \sqrt{\Delta (G)}\rceil +1$. 
We prove a more general result for the distinguishing index of graphs with  minimum degree at least two from which the conjecture follows.
Also  we present graphs $G$ for which $D'(G)\leq \lceil  \sqrt{\Delta }\rceil$.
\end{abstract}


\medskip
\noindent{\bf AMS Subj.\ Class.}: 05C25, 05C15

\section{Introduction }

Let $G=(V,E)$ be a simple connected graph. We use the standard graph notation (\cite{Sandi}). In particular, ${\rm Aut}(G)$ denotes the automorphism group of $G$.   For simple connected graph $G$, and $v \in V$,
the \textit{neighborhood} of a vertex $v$ is the set $N_G(v) = \{u \in V(G) : uv \in   E(G)\}$. The \textit{degree of a vertex} $v$ in a graph $G$, denoted by ${\rm deg}_G(v)$, is the number of edges of $G$ incident with $v$. In particular, ${\rm deg}_G(v)$ is the number of neighbours of $v$ in $G$.  We denote by $\delta(G)$ and $\Delta(G)$ the minimum and maximum degrees of the vertices of $G$. A graph $G$ is \textit{$k$-regular} if ${\rm deg}_G(v) = k$ for all $v \in V$. The \textit{diameter of a graph}  $G$ is the greatest distance between two vertices of $G$, and denoted by ${\rm diam}(G)$.

The \textit{distinguishing index} $D'(G)$ of a graph $G$ is the least number $d$ such that $G$ has an edge labeling  with $d$ labels that is preserved only by the identity automorphism of $G$.  The distinguishing edge labeling was first defined by  Kalinowski and Pi\'sniak \cite{R. Kalinowski and M. Pilsniak} for graphs (was inspired by the well-known distinguishing number $D(G)$ which was defined for general vertex labelings by Albertson and Collins \cite{Albert}). The distinguishing index of some examples of graphs was exhibited in \cite{R. Kalinowski and M. Pilsniak}. For 
 instance, $D'(P_n)=2$ for every $n\geq 3$, and 
 $ D'(C_n)=3$ for $n =3,4,5$,  $D'(C_n)=2$ for $n \geq 6$. Also, for complete graphs $K_n$, we have  $ D'(K_n)=3$ for $n =3,4,5$,  $D'(K_n)=2$ for $n \geq 6$.  They showed that if $G$ is a connected graph of order $n\geq 3$ and maximum degree $\Delta$,  then  $D'(G) \leq \Delta$, unless $G$ is $C_3, C_4$ or $C_5$.  It follows for connected graphs that $D'(G) \geq \Delta(G)$ if and only if $D'(G) =\Delta(G) + 1$ and $G$ is a cycle of length at most five. The equality $D'(G) =
\Delta(G)$ holds for all paths, for cycles of length at least 6, for $K_4$, $K_{3,3}$ and for
symmetric or bisymmetric trees. Also,  Pil\'sniak showed that $D'(G) < \Delta(G)$ for all
other connected graphs.
\begin{theorem}{\rm\cite{nord}}\label{updisindg}
Let $G$ be a connected graph that is neither a symmetric nor
an asymmetric tree. If the maximum degree of $G$ is at least 3, then $D'(G) \leq  \Delta(G) - 1$ unless $G$ is $K_4$ or $K_{3,3}$.
 \end{theorem}

  Pil\'sniak  put forward the following conjecture.
  
  \begin{conjecture} {\rm\cite{nord}}\label{konj}
  	If $G$ is a $2$-connected graph, then $D'(G)\leq 1+\lceil  \sqrt{\Delta (G)}\rceil$. 	
  \end{conjecture}

In this paper, we prove the following theorem which proves the conjecture. 
  \begin{theorem}
  	Let $G$ be a connected graph of maximum degree $\Delta$. If  the minimum degree $\delta \geq 2$, then $ D'(G)\leq \lceil  \sqrt{\Delta }\rceil +1$.
  \end{theorem}
  
   For our purposes, we consider graphs with specific construction that are from dutch-windmill graphs. Because of this,  in Section 2, we compute the distinguishing index of the dutch windmill graphs. In Section 3, we use the results to prove the main result.   In the last section we present graphs $G$ for which $D'(G)\leq \lceil  \sqrt{\Delta }\rceil$.

\section{Distinguishing index of dutch windmill graphs}

\begin{figure}
\begin{center}
\includegraphics[width=0.75\textwidth]{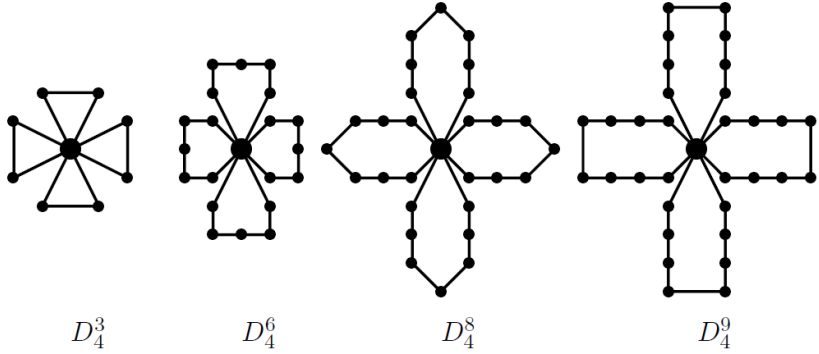}
	\caption{\label{fig0}  Examples of dutch windmill graphs.}
\end{center}
\end{figure}

To obtain the upper bound for the distinguishing index of connected graphs with minimum degree at least two, we characterize such graphs with minimum number of edges. For this characterization we need the concept of dutch windmill graphs. 
The \textit{dutch windmill graph}  $D_n^k$  is the graph obtained by taking $n$, ($n\geq 2$) copies of the cycle graph $C_k$, ($k\geq 3$) with a vertex in common (see Figure \ref{fig0}). If $k=3$, then we call $D_n^3$, a friendship graph. In the following theorem we compute the distinguishing number of dutch windmill graphs.  

\begin{figure}[ht]
\begin{center}
\includegraphics[width=0.9\textwidth]{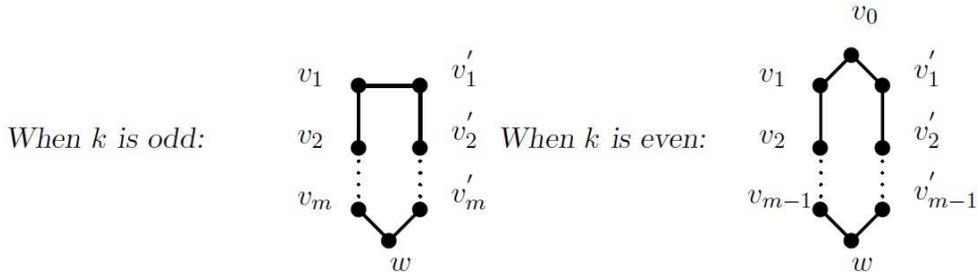}
	\caption{\label{fig00} {\small The considered polygon (or a cycle of size $k$) in the proof of Theorem \ref{thmbasic}}.}
\end{center}
\end{figure}

\begin{theorem}\label{thmbasic}
For every $n\geq 2$ and $k\geq 3$, 
$D(D_n^k)={\rm min}\{r:~\dfrac{r^{k-1}-r^{\lceil \dfrac{k-1}{2} \rceil}}{2}\geq n\}$.
\end{theorem}
\proof  
We consider two cases:

{\bf Case 1)} If $k$ is odd. There is a natural number $m$ such that $k=2m+1$. We can consider a blade of $D_n^k$ as Figure \ref{fig00}.
 Let $(x_1^{(i)},x_1'^{(i)},\ldots ,x_m^{(i)},x_m'^{(i)})$ be the label of vertices $(v_1,v_1',\ldots ,v_m,v_m')$ of the $i$th blade where $1\leq i \leq n$. Suppose that $L=\{(x_1^{(i)},x_1'^{(i)},\ldots ,x_m^{(i)}, x_m'^{(i)})\vert ~ 1\leq i \leq n , x_j^{(i)}, x_j'^{(i)}\in \mathbb{N}, 1\leq j \leq m \}$ is a labeling of the vertices of $D_n^k$ except its central vertex. In an $r$-distinguishing labeling we must have:
\begin{itemize}
\item[(i)] There exists $j\in \{1,\ldots ,m\}$ such that $x_j^{(i)}\neq x_j'^{(i)}$ for all $i\in \{1,\ldots , n\}$.
\item[(ii)] For $i_1\neq i_2$ we must have $(x_1^{(i_1)},x_1'^{(i_1)},\ldots ,x_m^{(i_1)},x_m'^{(i_1)})\neq (x_1^{(i_2)},x_1'^{(i_2)},\ldots ,x_m^{(i_2)},x_m'^{(i_2)})$  and  $(x_1^{(i_1)},x_1'^{(i_1)},\ldots ,x_m^{(i_1)},x_m'^{(i_1)})\neq (x_1'^{(i_2)},x_1^{(i_2)},\ldots ,x_m'^{(i_2)},x_m^{(i_2)})$. 
\end{itemize}
There are $\dfrac{r^{2m}-r^m}{2}$ possible $(2m)$-arrays of labels using $r$ labels satisfying (i) and (ii), and so  $D(D_n^k)={\rm min}\{r:~\dfrac{r^{2m}-r^m}{2}\geq n\}$.

{\bf Case 2)} If $k$ is even. There is a natural number $m$ such that $k=2m$. We can consider a blade of $D_n^k$ as Figure \ref{fig00}. Let $(x_0^{(i)}x_1^{(i)}, x_1'^{(i)}, \ldots , x_{m-1}^{(i)}, x_{m-1}'^{(i)})$ be the label of vertices $(v_0,v_1,v_1',\ldots ,v_{m-1},v_{m-1}')$ of  $i$th blade where $1\leq i \leq n$. Suppose that $L=\{(x_0^{(i)},x_1^{(i)},x_1'^{(i)},\ldots ,x_{m-1}^{(i)}, x_{m-1}'^{(i)})\vert ~ 1\leq i \leq n , x_0^{(i)},x_j^{(i)},x_j'^{(i)}\in \mathbb{N}, 1\leq j \leq m-1 \}$ is a labeling of the vertices of $D_n^k$ except its central vertex. In an $r$-distinguishing labeling we must have:
\begin{itemize}
\item[(i)] There exists $j\in \{1,\ldots ,m-1\}$ such that $x_j^{(i)}\neq x_j'^{(i)}$ for all $i\in \{1,\ldots , n\}$.
\item[(ii)] For $i_1\neq i_2$ we must have 
\begin{align*}
&(x_0^{(i_1)},x_1^{(i_1)},x_1'^{(i_1)},\ldots ,x_{m-1}^{(i_1)},x_{m-1}'^{(i_1)})\neq (x_0^{(i_2)},x_1^{(i_2)},x_1'^{(i_2)},\ldots ,x_{m-1}^{(i_2)},x_{m-1}'^{(i_2)}),\\ &(x_0^{(i_1)},x_1^{(i_1)},x_1'^{(i_1)},\ldots ,x_{m-1}^{(i_1)},x_{m-1}'^{(i_1)})\neq (x_0^{(i_2)},x_1'^{(i_2)},x_1^{(i_2)},\ldots ,x_{m-1}'^{(i_2)},x_{m-1}^{(i_2)}).
\end{align*}
\end{itemize}
There are $\dfrac{r^{2m-1}-r^m}{2}$ possible $(2m-1)$-arrays of labels using $r$ labels satisfying (i) and (ii) ($r$ choices for $x_0$ and $\dfrac{r^{2(m-1)}-r^{m-1}}{2}$ choices for $x_1^{(i_1)},x_1'^{(i_1)},\ldots ,x_{m-1}^{(i_1)},x_{m-1}'^{(i_1)}$). Therefore  $D(D_n^k)={\rm min}\{r:~\dfrac{r^{2m-1}-r^m}{2}\geq n\}$.\qed

The following theorem gives the distinguishing index of  $D_n^k$.  
\begin{theorem}\label{indexwindmil}
For any $n\geq 2$ and $k\geq 3$,  $D'(D_n^k)={\rm min}\{r:~\dfrac{r^{k}-r^{\lceil \dfrac{k}{2} \rceil}}{2}\geq n\}$.
\end{theorem}
\proof
Since the effect of every automorphism of $D_n^{k+1}$ on its non-central vertices is exactly the same as the effect of an automorphism of $D_n^k$ on its edges and vice versa, so if we consider the non-central vertices of $D_n^{k+1}$ as the edges of $D_n^k$, then we have $D'(D_n^k)=D(D_n^{k+1})$. Therefore the result follows from Theorem \ref{thmbasic}.\qed

\section{Proof of conjecture}

In this section, we shall prove  Conjecture \ref{konj}. To do this, first we state some preliminaries.  
By the result obtained  by Fisher and Isaak \cite{fish} and independently by Imrich, Jerebic and Klav\v zar \cite{W Imrich} the distinguishing index of complete bipartite graphs is as  follows:
\begin{theorem}{\rm \cite{fish, W Imrich}}\label{indcombipar}
	Let $p, q, d$ be integers such that $d \geq 2$ and $(d-1)^p <
	q \leq d^p$ . Then
	\begin{equation*}
	D'(K_{p,q}) =\left\{
	\begin{array}{ll}
	d & \text{if}~~ q \leq d^p - \lceil {\rm log}_d p\rceil - 1,\\
	d + 1 & \text{if}~~ q \geq d^p - \lceil {\rm log}_d p\rceil + 1.
	\end{array}\right.
	\end{equation*}
	
	If $q = d^p -\lceil {\rm log}_d p\rceil$ then the distinguishing index $D'(K_{p,q})$ is either $d$ or $d+1$ and can be computed recursively in $O({\rm log} ∗(q))$ time.
\end{theorem}

\begin{corollary}{\rm \cite{nord}}\label{uppcompbipart}
	If $p \leq q$, then $D'(K_{p,q}) \leq \lceil \sqrt[p]{q} \rceil + 1$.
\end{corollary}

  Also we need the following result:
 \begin{theorem}{\rm\cite{nord}}\label{hamiltopath}
 If $G$ is a  graph of order $n \geq 7$ such that $G$ has a Hamiltonian path, then $D'(G) \leq 2$.
  \end{theorem}
 
 Now, we state and prove the main theorem of this paper. 

\begin{theorem}\label{uppinddelta2}
Let $G$ be a connected   graph with maximum degree $\Delta$. If $\delta \geq 2$ then $D'(G) \leq 1+ \lceil \sqrt{\Delta}\rceil$.
\end{theorem}
\proof  If $\Delta \leq 5$, then the result follows from Theorem \ref{updisindg}. So, we suppose that $\Delta \geq 6$. 
 Let $v$ be a vertex of $G$ with the maximum degree $\Delta$.  By Theorem \ref{indexwindmil},  we can label the edges of the  dutch windmill  graph attached to $G$ at vertex $v$ (a subgraph $H$ is attached  to graph $G$, if   it has only
one vertex in common with graph $G$)  for which $v$ is the central point of the dutch windmill  graph, with at most $ \lceil \sqrt{\Delta}\rceil$ labels from label set $\{0, 1, \ldots ,  \lceil \sqrt{\Delta}\rceil\}$, distinguishingly.   If there exists triangle attached to  $G$ at $v$, then we label the two its incident edges to $v$ with $0$ and $1$, and another edges of the  triangle with label $2$. 

Let $N^{(1)}(v) =\{v_1, \ldots , v_{|N^{(1)}(v) |}\}$ be the vertices of $G$ at distance one from $v$, except the vertices of  dutch windmill  or triangle attached to graph at $v$. We continue the  labeling by the following steps:

\medskip

{ \bf Step 1)} Since $|N^{(1)}(v)| \leq \Delta$, so for $0 \leq i \leq \lceil \sqrt{\Delta}\rceil$ and $1 \leq j \leq \lceil \sqrt{\Delta}\rceil-1$,  we  label the edges $vv_{i \lceil \sqrt{\Delta}\rceil +j}$ with label $i$,  and we do not use the label $0$ any more. With respect to the number of incident edges to $v$ with label 0, we conclude that the vertex $v$ is fixed under each automorphism of $G$ preserving the labeling. Also, since the dutch windmill  or the triangle graph attached  to $G$ at $v$ has been labeled distinguishingly, so the vertices of attached  graph  are fixed under each automorphism of $G$ preserving the labeling. Hence, every automorphism of $G$ preserving the labeling must map the set of vertices of $G$ at distance $i$ from $v$ to itself setwise,  for any $1 \leq i \leq {\rm diam}(G)$. We denote the set of vertices of $G$ at distance $i$ from $v$, by $N^{(i)}(v)$ 
for $2 \leq i \leq {\rm diam}(G)$.  If for any $i \geq 2$, $N^{(i)}(v) = \emptyset$,  then $G$ has a Hamiltonian path, and since $\Delta \geq 6$, so the order of $G$ is at least $7$, and hence $D'(G) \leq 2$ by Theorem \ref{hamiltopath}. Thus we suppose that $N^{(i)}(v) \neq  \emptyset$, for some $i \geq 2$. 

 Now we partition the vertices of $N^{(1)}(v)$ to two sets $M_1^{(1)}$ and $M_2^{(1)}$ as follows:
 \begin{equation*}
 M_1^{(1)} =\{x\in N^{(1)}(v):~ N(x) \subseteq N(v)\},~~M_2^{(1)}=\{x\in N^{(1)}(v):~ N(x) \nsubseteq N(v)\}.
 \end{equation*}
 The sets $M_1^{(1)}$ and $M_2^{(1)}$ are mapped to $M_1^{(1)}$ and $M_2^{(1)}$, respectively, setwise, under each automorphism of $G$ preserving the labeling. For  $0 \leq i \leq \lceil \sqrt{\Delta}\rceil$, we set $L_i =\{v_{i \lceil \sqrt{\Delta}\rceil +j}~:~ 1 \leq j \leq \lceil \sqrt{\Delta}\rceil -1\}$. By this notation, we get that for $0 \leq i \leq \lceil \sqrt{\Delta}\rceil$, the set $L_i$ is mapped to $L_i$ setwise,  under each automorphism of $G$ preserving the labeling. Let the sets $M_{1i}^{(1)}$ and $M_{2i}^{(1)}$ for $0 \leq i \leq \lceil \sqrt{\Delta}\rceil$ are as follows:
\begin{equation*}
M_{1i}^{(1)}= M_1^{(1)} \cap L_i,~~ M_{2i}^{(1)}= M_2^{(1)} \cap L_i.
\end{equation*}

It is clear that the sets $M_{1i}^{(1)}$ and $M_{2i}^{(1)}$  are mapped to $M_{1i}^{(1)}$ and $M_{2i}^{(1)}$, respectively, setwise,  under each automorphism of $G$ preserving the labeling. Since  for any $0 \leq i \leq \lceil \sqrt{\Delta}\rceil$, we have $|M_{1i}^{(1)}|\leq \lceil \sqrt{\Delta}\rceil-1$, so we can label all incident edges to each element of $M_{1i}^{(1)}$ with labels $\{1,2, \ldots , \lceil \sqrt{\Delta}\rceil \}$, such that for any two vertices of $M_{1i}^{(1)}$, say $x$ and $y$, there exists a label $k$, $1\leq k \leq \lceil \sqrt{\Delta}\rceil$, such that the number of label $k$ for the incident edges to vertex $x$ is different from the number of label $k$ for the incident edges to vertex  $y$. Hence, it can be deduce that each vertex of  $M_{1i}^{(1)}$ is fixed  under each automorphism of $G$ preserving the labeling, where $0\leq i \leq  \lceil \sqrt{\Delta}\rceil$. Thus every vertices of $M_1^{(1)}$ is fixed  under each automorphism of $G$ preserving the labeling. In sequel, we want to label the edges incident to vertices of $M_2^{(1)}$ such that $M_2^{(1)}$ is fixed  under each automorphism of $G$ preserving the labeling, pointwise. For this purpose, we partition the vertices of $M_{2i}^{(1)}$ to the sets $M_{{2i}_j}^{(1)}$, ($1 \leq j \leq \Delta -1$) as follows:
\begin{equation*}
M_{{2i}_j}^{(1)} =\{x\in M_{2i}^{(1)}~:~ |N(x) \cap  N^{(2)}(v)| = j\}.
\end{equation*}

Since the set $N^{(i)}(v)$, for any $i$, is mapped to itself, it can be concluded that $M_{{2i}_j}^{(1)}$  is mapped to itself, setwise, under each automorphism of $G$ preserving the labeling, for any $i$ and $j$.  Let $M_{{2i}_j}^{(1)} =\{x_{i_{j1}},\ldots , x_{i_{js_j}}\}$. It is clear that $|M_{{2i}_j}^{(1)}| \leq |M_{2i}^{(1)}| \leq \lceil \sqrt{\Delta}\rceil-1$. Let $x_{i_{jk}}\in M_{{2i}_j}^{(1)}$, and $N(x_{i_{jk}}) \cap N^{(2)}(v) =\{x'_{i_{jk1}}, x'_{i_{jk2}}, \ldots , x'_{i_{jkj}}\}$. We assign to the  $j$-ary $(x_{i_{jk}}x'_{i_{jk1}}, \ldots , x_{i_{jk}}x'_{i_{jkj}})$ of edges, a $j$-ary of labels such that for every $x_{i_{jk}}$ and $x_{i_{jk'}}$, $1\leq k,k'\leq s_j$, there exists a label $l$ in their corresponding $j$-arys of labels for which the  number of label $l$ in the corresponding $j$-arys of $x_{i_{jk}}$ and $x_{i_{jk'}}$ is distinct. For constructing $| M_{{2i}_j}^{(1)}|$ numbers of such $j$-arys we need, ${\rm min}\{r:~ {j+r-1 \choose r-1}\geq | M_{{2i}_j}^{(1)}| \}$ distinct labels. Since for any $1 \leq j \leq \Delta -1$, we have
\begin{equation*}
{\rm min}\left\{r:~ {j+r-1 \choose r-1}\geq | M_{{2i}_j}^{(1)}| \right\} \leq {\rm min}\left\{r:~ {j+r-1 \choose r-1}\geq \lceil \sqrt{\Delta}\rceil -1 \right\} \leq \lceil \sqrt{\Delta}\rceil,
\end{equation*}
so we need at most $\lceil \sqrt{\Delta}\rceil$ distinct labels from label set $\{1, 2, \ldots , \lceil \sqrt{\Delta}\rceil\}$ for constructing such $j$-arys. For instance, let $j=1$, and $M_{{2i}_1}^{(1)}=\{x_{i_{11}}, \ldots , x_{i_{1s_1}}\}$. By our method, we label the edge $x_{i_{11}}x'_{i_{1k1}}$ with label $k$ for $1 \leq k \leq s_1$ where $s_1 \leq  \lceil \sqrt{\Delta}\rceil$. Hence, the vertices of $M_{{2i}_j}^{(1)}$, for any $1 \leq j \leq \Delta -1$, are fixed under each automorphism of $G$ preserving  the labeling. Therefore, the vertices of $M_{2i}^{(1)}$ for any $0 \leq i \leq \lceil \sqrt{\Delta}\rceil$, and so the vertices of $M_2^{(1)}$ are fixed under each automorphism of $G$ preserving  the labeling. Now, we can get that all vertices of $N^{(1)}(v)$ are fixed. If there exist unlabeled edges of $G$ with the two endpoints in $N^{(1)}(v)$, then we assign them an arbitrary label, say 1.

\medskip
{ \bf Step 2)} Now we consider $N^{(2)}(v)$. We partition this set such that the vertices of $N^{(2)}(v)$ with the same neighbours in $M_2^{(1)}$, lie in a set. In other words, we can write $N^{(2)}(v) = \bigcup_i A_i$, such that $A_i$ contains that elements of $N^{(2)}(v)$ having the same neighbours in $M_2^{(1)}$, for any $i$. Since all vertices in  $M_2^{(1)}$  are fixed, so  the set $A_i$ is mapped to $A_i$ setwise, under each automorphism of $G$ preserving the labeling. Let $A_i = \{w_{i1}, \ldots , w_{it_i}\}$, and we have 
\begin{equation*}
N(w_{i1}) \cap M_2^{(1)} = \cdots = N(w_{it_i}) \cap M_2^{(1)} = \{v_{i1}, \ldots , v_{ip_i}\}.
\end{equation*}

We consider two following cases:

Case 1) If for every $w_{ij}$ and $w_{ij'}$ in $A_i$, where $1\leq j,j' \leq t_i$, there exists a $k$, $1\leq k \leq p_i$, for which the label of edge $w_{ij}v_{ik}$ is different from label of edge $w_{ij'}v_{ik}$, then all vertices of $G$ in $A_i$ are fixed under each automorphism  of $G$ preserving the labeling.

Case 2) If there exist  $w_{ij}$ and $w_{ij'}$ in $A_i$, where $1\leq j,j' \leq t_i$, such that for every  $k$, $1\leq k \leq p_i$,  the label of edgeس $w_{ij}v_{ik}$ and $w_{ij'}v_{ik}$ are the same, then we can make a labeling  such that the  vertices  in $A_i$ have the same property   as Case 1, and so are fixed under each automorphism  of $G$ preserving the labeling, by using at least one of the following actions:
\begin{itemize}
\item By commutating the coordinates of $j$-ary of labels assigned to the incident edges to $v_{ik}$ with an end point in $N^{(2)}(v)$. 
\item By using a new $j$-ary of labels, with labels $\{1, 2, \ldots , \lceil \sqrt{\Delta}\rceil\}$,  for incident edges to $v_{ik}$ with an end point in $N^{(2)}(v)$, such that (by notations in Step 1)
 for every $x_{i'_{jk'}}$ and $x_{i'_{jk''}}$, $1\leq k',k''\leq s_j$, there exists a label $l$ in their corresponding $j$-arys of labels with different number of label $l$ in their coordinates, where $1 \leq i' \leq \lceil \sqrt{\Delta} \rceil$. 
\item  By labeling the unlabeled edges of $G$ with the two end points in $N^{(2)}(v)$ which are incident  to the vertices  in $A_i$. 
\item By labeling the unlabeled edges of $G$   which are incident  to the vertices  in $A_i$, and another their endpoint is $N^{(3)}(v)$.
\item By labeling the unlabeled edges of $G$ with the two end points in $N^{(3)}(v)$ for which the end points in $N^{(3)}(v)$ are adjacent  to some of  vertices  in $A_i$.
\end{itemize}

Using at least one of above actions, it can be seen that every two vertices $w_{ij}$ and $w_{ij'}$ in $A_i$ have the property as Case 1. Thus we conclude that all vertices in $A_i$, for any $i$, and so all vertices in $N^{(2)}(v)$, are fixed under each automorphism  of $G$ preserving the labeling.  If there exist unlabeled edges of $G$ with the two endpoints in $N^{(2)}(v)$, then we assign them an arbitrary label, say 1.

By continuing  this method,  in the next step we  partition $N^{(3)}(v)$ exactly by the same method as partition of  $N^{(2)}(v)$ to the sets $A_i$'s  in  Step 2, and so we can make a labeling such that $N^{(i)}(v)$ is fixed pointwise, under each automorphism  of $G$ preserving the labeling, for any $3 \leq i \leq {\rm diam}(G)$.\qed

\medskip

For a $2$-connected planar graph $G$, the distinguishing index may attain
$1+\lceil \sqrt{\Delta (G)} \rceil$. For example, consider the complete bipartite graph $K_{2,q}$ with $q = r^2$, where $r$ is a positive integer $r$. By   Theorem \ref{indcombipar},  $D'(K_{2,q}) = r + 1$. 


\section{Graphs with $D'(G) \leq \lceil \sqrt{\Delta}\rceil$}

 In this section,  we present graphs $G$ with specific construction such that   $D'(G)\leq \lceil \sqrt{\Delta}\rceil$. To do this we state the following definition. 
\begin{definition}\label{minimally}
Let $G$ be a connected graph with $\delta(G) \geq 2$. The graph $G$ is called  a {\rm $\delta$-minimally graph}, if the minimum degree of each spanning subgraph of $G$ is less than $\delta(G)$. 
\end{definition}

It can be concluded from Definition \ref{minimally} that if $e$ is an edge of a connected  $\delta$-minimally graph with end points $u$ and $w$, then without loss of generality we can assume that ${\rm deg}_G u= \delta$ and ${\rm deg}_G w\geq \delta$. In fact the distance between the two vertices of degree greater than $\delta$ is at least two. 

The simplest connected  $2$-minimally graphs are cycles $C_n$ and complete bipartite graphs $K_{2,n}$. Now, we explain more on the structure of a $2$-minimally graph.
Let to call  a path in the graph a {\it simple path}, if  all its internal vertices have degree two.  Let $G$ be a connected  $2$-minimally graph. 
\begin{itemize}
\item If the degree of all vertices of $G$ is two, then $G$ is a cycle graph. 

\item If there exist a vertex $v$ of $G$ with degree at least three. We consider  two following cases:

Case 1)   If $v$ is the only vertex of $G$ with degree greater than two, then $G$ is a graph which is made by identifying the central points of some dutch windmill graphs $D_{n_i}^{p_i}$ where $p_i\geq 3$, and hence $\Delta(G) = 2\sum_{i\in I} n_i$ where $I$ is a set of indices. In this case we denote $G$ by ${\rm Wind}(v)$ (for instance, see Figure \ref{figwind}).

\begin{figure}
\begin{center}
\includegraphics[width=0.8\textwidth]{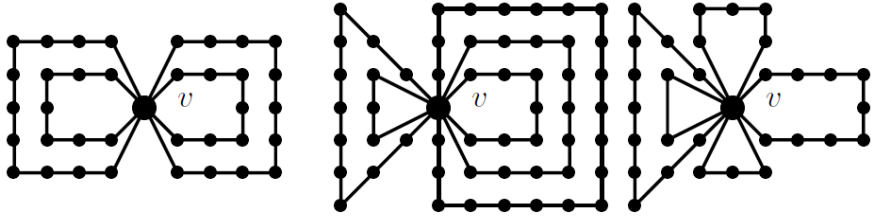}
	\caption{\label{figwind}  Examples of ${\rm Wind}(v)$.}
\end{center}
\end{figure}

Case 2) If $G$ has other vertex $w$ of degree greater than two, then there exists  at least  a simple path between $v$ and $w$ of length	greater than one. Since $G$ is a connected graph, so if there exists no such simple path, then there exists a vertex  of degree at least three on each path between $v$ and $w$. Hence we can obtain a vertex $u$ of $G$ with degree greater than two such that there exists at least a simple  path between $v$ and $u$ of length  greater than one (see Figure \ref{fig1}).
\end{itemize}

\begin{figure}[ht]
	\hspace{1cm}
	\begin{minipage}{6.2cm}
		\includegraphics[width=1.0\textwidth]{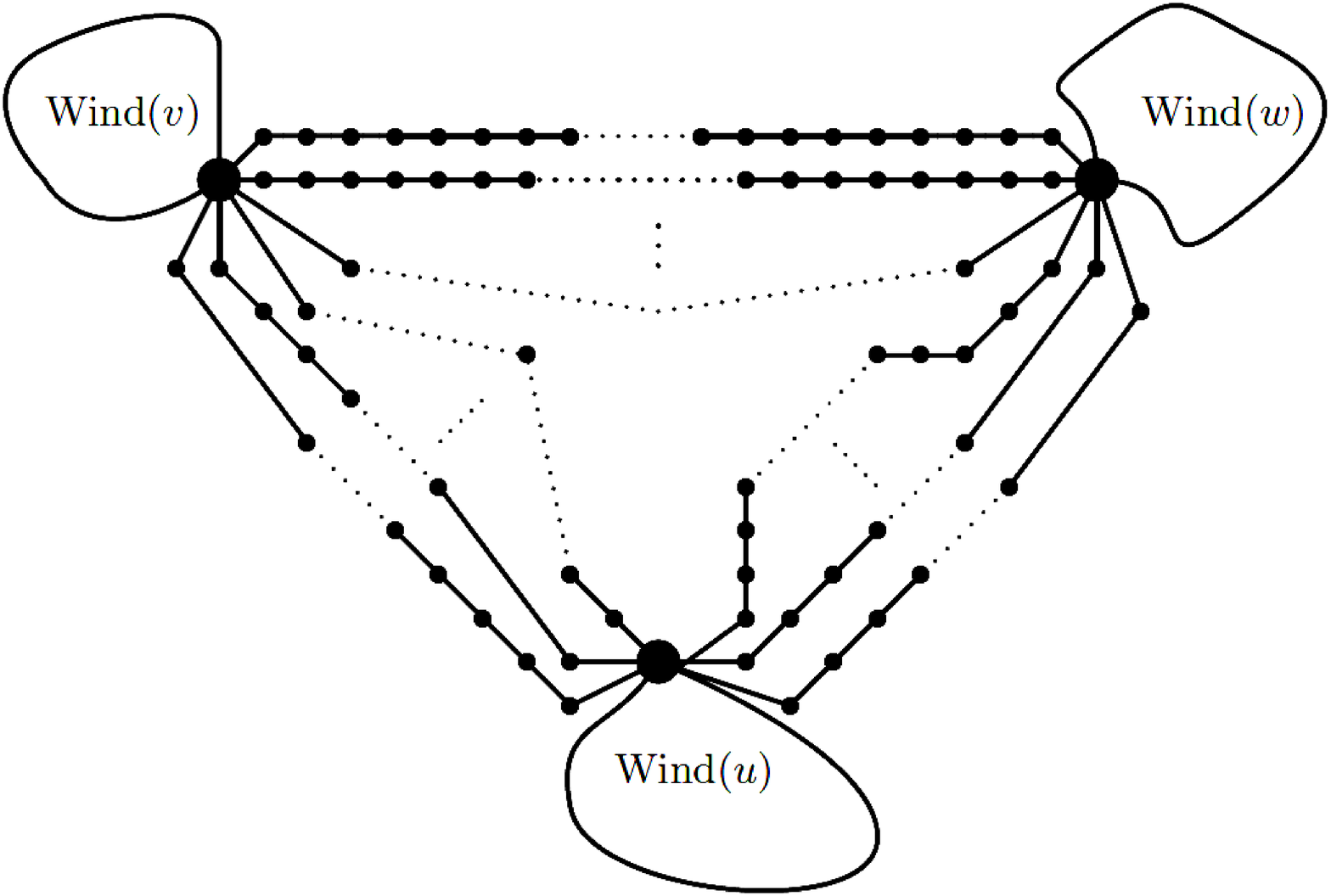}
	\end{minipage}
	\begin{minipage}{6.2cm}
		\includegraphics[width=1.0\textwidth]{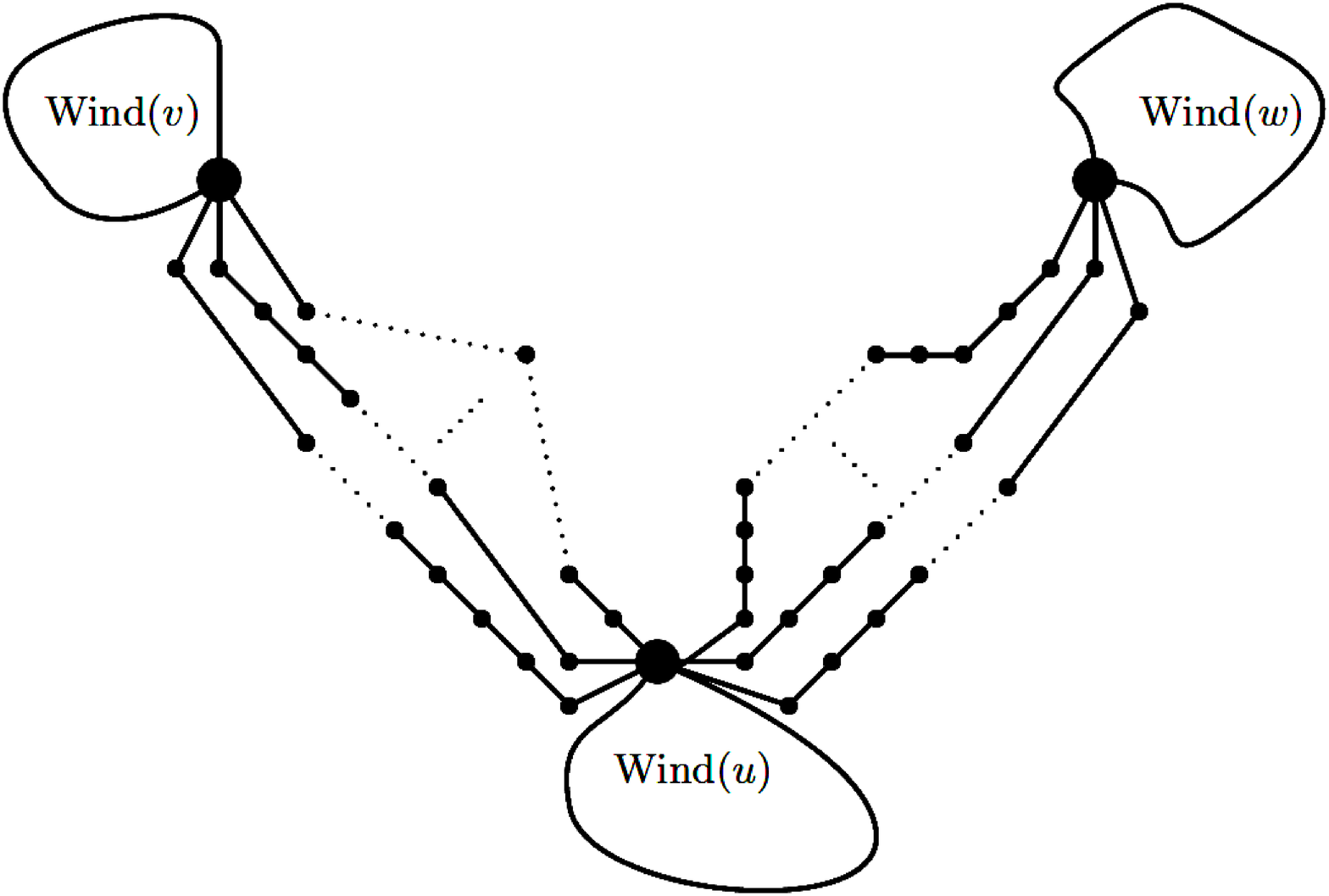}
	\end{minipage}
	\caption{\label{fig1} The state of vertices of degree greater than two in a connected $2$-minimally graph.}
\end{figure}

Now we characterize graphs $G$ with $D'(G) \leq \lceil \sqrt{\Delta} \rceil$. 
\begin{theorem}\label{indexminimallygraphs}
Let $G$ be a connected  $2$-minimally graph with maximum degree $\Delta$. If $G$ is not a cycle $C_3, C_4, C_5$ or a complete bipartite graph $K_{2, r^2}$ for some integer $r$, then $D'(G) \leq \lceil \sqrt{\Delta} \rceil$. 
\end{theorem}
\proof
If $\Delta = 2$, then $G$ is a cycle. It is known that  the distinguishing index of cycle graph of order at least $6$ is two. Hence, we suppose that  $G$ is not a cycle, so $\Delta \geq 3$. 
Let $v$ be a vertex of $G$ of maximum degree $\Delta$.  Suppose that $V'=\{v_1,v_2, \ldots , v_k\}$ are all vertices of $G$ which are of degree at least three such that there exists at least a simple path between $v$ and $v_i$, for any $1 \leq i \leq k$ (it is possible that $V' =\emptyset$).  Let there exist $n_{ij}$ disjoint simple paths of length $j$ between $v$ and $v_i$,  for any $1 \leq i \leq k$ and $2 \leq j \leq {\rm diam}(G)$ where $n_{ij}$ is a non-negative integer and $\sum_{j=2}^{{\rm diam}(G)}n_{ij} > 0$.  We can label these $n_{ij}$ simple paths of length $j$ with at most $\lceil \sqrt{\Delta} \rceil$ labels, by using $n_{ij}$ numbers of $j$-arys such that the coordinates of each $j$-ary are in the set $\{1, 2, \ldots , \lceil \sqrt{\Delta} \rceil\}$,  for any $1\leq i \leq k$ and $2 \leq j \leq {\rm diam}(G)$, and for every two paths of length $j$, say $P_1$ and $P_2$, there exists a label $l$, $1 \leq l \leq \lceil \sqrt{\Delta} \rceil$, such that the number of label $l$ in $j$-arys related to $P_1$ and $P_2$ are distinct. Let $P$ be a simple path between $v$ and $v_i$ for some $1\leq i \leq k$, such that the label of edge of $P$ which is incident to $v$, is different from the label of edge of $P$ which is incident to $v_i$. We do not use of labeling of the simple path $P$, for any other simple path (with the same length) between any two vertices of degree greater than two. Since $G$ is not a complete bipartite graph $K_{2,r^2}$ for some integer $r$, so we can label these paths distinguishingly with at most $\lceil  \sqrt{\Delta} \rceil$ labels.  Now,   we label the induced subgraph ${\rm Wind} (x)$, for any vertex $x$ of degree greater than two, if there exists,  with at most $\lceil \sqrt{\Delta} \rceil$ labels distinguishingly by Theorem \ref{indexwindmil}, such that the distinguishing  labeling of ${\rm Wind} (v)$ is nonisomorphic to the remaining  distinguishing  labeling of ${\rm Wind} (x)$, where $x\in V(G)-\{v\}$. Thus any automorphism of $G$ preserving this labeling should be fixed $v, v_1,\ldots , v_k$ and all vertices of degree two on the simple paths  between $v$ and $v_i$ for any $1\leq i \leq k$. Since for any $1\leq i,j \leq k$ where $i \neq j$, the vertices $v_i$ and $v_j$ are fixed, so all  the simple paths between $v_i$ and $v_j$, if there exist, are mapped to each other  under each automorphism of $G$ preserving this labeling. Hence we can label all  edges of these simple paths with at most $\sqrt{\Delta}$ labels, by assigning distinct ordered arys of labels of length of  the simple paths between $v_i$ and $v_j$ such that  all vertices of these paths are fixed under each automorphism of $G$ preserving this labeling.

For any $1\leq i \leq k$, we consider $v_i$, and suppose that $v_{i1}, \ldots , v_{ik_i}$ are all vertices of $V(G) \setminus \{v_1,\ldots , v_k\}$ with degree at least three such that there exists at least a simple path between $v_i$ and $v_{ij}$ for any  $1\leq j\leq k_i$. Now we do the same method as labeling of simple paths between  $v$ and $\{v_1, \ldots , v_k\}$, for all simple paths  between  $v_i$ and $\{v_{i1}, \ldots , v_{ik_i}\}$ with at most $\lceil \sqrt{\Delta}  \rceil$ labels. Also, we do the same method as labeling of simple paths between  $v_i$ and $v_j$, for all simple paths  between  $v_{ip}$ and $v_{iq}$ with at most $\lceil \sqrt{\Delta}\rceil$ labels, where $1 \leq p,q\leq k_i$.  Note that we do not use  labeling of $P$ for any simple path with the same length as $P$ between $v_i$ and $\{v_{i1}, \ldots , v_{ik_i}\}$. Thus the vertices $\{v_i, v_{i1}, \ldots , v_{ik_i}\}$ and all  vertices of the simple paths between them are fixed under each automorphism of $G$ preserving this labeling. After the finite number of steps we can obtain a distinguishing edge labeling of $G$ with at most $\lceil \sqrt{\Delta} \rceil$ labels.\qed

\section{Conclusion} 

We gave an upper bound for the distinguishing index of graphs $G$ with minimum degree at least two. This result proves a conjecture by Pil\'sniak (2015). We also studied graphs $G$ with  $D'(G) \leq \lceil \sqrt{\Delta} \rceil$. 
  We think that the following conjecture is true, but    until now all attempts to prove this  failed. So, we  end  this paper by proposing  the following conjecture. 

\begin{conjecture} 
	Let $G$ be a connected graph with maximum and minimum degree $\Delta$ and $\delta$, respectively.
\begin{itemize}
\item[(i)] If $G$ is a $\delta$-minimally graph with $\delta(G)\geq 3$ such that $G$ is not a complete bipartite  or $\delta$-regular graph, then $D'(G)  \leq \lceil  \sqrt[\delta(G)]{\Delta (G)}\rceil$.
\item[(ii)] If $G$ is a connected graph with $\delta(G)\geq 3$, then $D'(G)  \leq 1+ \lceil  \sqrt[\delta(G)]{\Delta (G)}\rceil$.
\end{itemize}
\end{conjecture}


\begin{thebibliography}{99}

	\bibitem{Albert} M.O. Albertson and K.L. Collins, {\it Symmetry breaking in graphs}, Electron. J. Combin. 3 (1996), \#R18.


\bibitem{soltani2} S. Alikhani and S. Soltani, {\it Distinguishing number and distinguishing index of certain graphs}, Filomat, to appear. http://arxiv.org/abs/1602.03302.


\bibitem{fish} M. J. Fisher and G. Isaak, {\it Distinguishing colorings of Cartesian products of complete graphs}, Discrete Math. 308 (11) (2008), 2240-2246. 
\bibitem{Sandi} R. Hammack, W. Imrich and S. Kla\v{v}zar, {\it Handbook of product graphs (second edition)}, Taylor \& Francis group, 2011.  



\bibitem{W Imrich} W. Imrich, J. Jerebic and S. Klav\v zar, {\it The distinguishing number of
	Cartesian products of complete graphs}, European J. Combin. 29 (2008),
922-929.
\bibitem{R. Kalinowski and M. Pilsniak} R. Kalinowski and M. Pil\'sniak, {\it Distinguishing graphs by edge colourings}, 
European J. Combin. 45 (2015), 124-131.



\bibitem{nord}  M. Pil\'sniak, {\it Nordhaus-Gaddum bounds for the distinguishing index}, (2015) Available at \texttt{www.ii.uj.edu.pl/preMD/}. 
\end{thebibliography}
\end{document}